\documentclass[11pt,reqno]{amsart}

\numberwithin{equation}{subsection}
\usepackage{framed} % for leftbar
\usepackage{tikz}
\usepackage{etex}
\usetikzlibrary{shapes,calc,arrows,through,intersections}

\tikzset{edge from parent/.append style={->}}
\usepackage{tikz-qtree}
\usepackage[margin=1in]{geometry}

\usepackage{stmaryrd}

\usepackage{amsmath,amscd,amssymb,amsfonts,latexsym,wasysym, mathrsfs, mathtools,hhline,color, bm}
\usepackage[all, cmtip]{xy}

\usepackage{url}

\definecolor{hot}{RGB}{65,105,225}

\usepackage[pagebackref=true,colorlinks=true, linkcolor=hot ,  citecolor=hot, urlcolor=hot]{hyperref}%make all the references and links clickable

\usepackage{ textcomp }
\usepackage{ tipa }

\usepackage{graphicx,enumerate}
\usepackage{subcaption}

\usepackage{enumitem}
\usepackage[normalem]{ulem}

\newcommand{\bfe}{\mathbf{e}}

\newcommand{\calV}{\mathcal{V}}

\newcommand{\ndot}{{n}_\bullet}

\theoremstyle{plain}
\newtheorem{theorem}{Theorem}[section]

\newtheorem{lm}[theorem]{Lemma}

\newtheorem{thrm}[theorem]{Theorem}

\theoremstyle{definition}
\newtheorem{alg}[theorem]{Algorithm}

\newtheorem{remark}[theorem]{Remark}

\newtheorem{ex}[theorem]{Example}
\newtheorem*{ex*}{Example}

\def\be{\begin{equation}}
\def\ee{\end{equation}}

\def\bt{\begin{thrm}}
\def\et{\end{thrm}}

\def\bc{\begin{cor}}
\def\ec{\end{cor}}

\def\br{\begin{rmk}}
\def\er{\end{rmk}}

\def\bp{\begin{prop}}
\def\ep{\end{prop}}

\def\bl{\begin{lm}}
\def\el{\end{lm}}

\def\bex{\begin{ex}}
\def\eex{\end{ex}}

\def\bd{\begin{defn}}
\def\ed{\end{defn}}

\newcommand\sV{{\mathcal V}}

\def\bfx{\mathbf{x}}

\newcommand\pp{{\mathbb{P}}}

\newcommand{\blue}[1]{{\color{blue}#1}}
\newcommand{\magenta}[1]{{\color{magenta}#1}}

\usepackage{algorithm}          %  float wrapper for algorithms.
\usepackage{algpseudocode}      % layout for algorithmicx

\algblock{Input}{EndInput}
\algnotext{EndInput}
\algblock{Output}{EndOutput}
\algnotext{EndOutput}
\algblock{Action}{EndAction}
\algnotext{EndAction}
\algblock{Compute}{EndCompute}
\algnotext{EndCompute}
\algblock{Initialize}{EndInitialize}
\algnotext{EndInitialize}
\algblock{Procedure}{EndProcedure}
\algnotext{EndProcedure}
\algblockdefx[Foreach]{Foreach}{EndForeach}[1]{\textbf{for each} #1 \textbf{do}}{\textbf{end for}}
\begin{document}

\title{MultiRegeneration for polynomial system solving}

%% Authors by alphabetical order. 

%%%%%%%%%%%%%%%%%%%%%%%%%%%%%%%%%%%%%%%%%%%%%%%%%%%%%%%%%%%%%%%%%%%%%%%%%%%% 
\author[C. Crowley]{Colin Crowley}
\address{Department of Mathematics\\
         University of Wisconsin---Madison\\
         Madison, WI 53706\\         
         USA}
\email{cwcrowley@wisc.edu}
\urladdr{{http://www.math.wisc.edu/~crowley/}}
%%%%%%%%%%%%%%%%%%%%%%%%%%%%%%%%%%%%%%%%%%%%%%%%%%%%%%%%%%%%%%%%%%%%%%%%%%%% 
\author[J.~I.~Rodriguez]{Jose Israel Rodriguez}
\address{Department of Mathematics\\
         University of Wisconsin---Madison\\
         Madison, WI 53706\\         
         USA}
\email{Jose@math.wisc.edu}
\urladdr{{http://www.math.wisc.edu/~jose/}}
%%%%%%%%%%%%%%%%%%%%%%%%%%%%%%%%%%%%%%%%%%%%%%%%%%%%%%%%%%%%%%%%%%%%%%%%%%%% 
\author[J.~Weiker]{Jacob Weiker}
\address{Department of Mathematics\\
         University of Wisconsin---Madison\\
         Madison, WI 53706\\         
         USA}
%\email{???}
%\urladdr{{???}}
%%%%%%%%%%%%%%%%%%%%%%%%%%%%%%%%%%%%%%%%%%%%%%%%%%%%%%%%%%%%%%%%%%%%%%%%%%%% 
\author[J.~Zoromski]{Jacob Zoromski}
\address{Department of Mathematics\\
		 University of Notre Dame\\
		 Notre Dame, IN 46556\\
                 USA}
\email{jzoromsk@nd.edu}
%\urladdr{{???}}
%%%%%%%%%%%%%%%%%%%%%%%%%%%%%%%%%%%%%%%%%%%%%%%%%%%%%%%%%%%%%%%%%%%%%%%%%%%%%%%%%

%%%%%%%%%%%%%%%%%%%%%%%%%%%%%%%%%%%%%%%%%%%%%%%%%%%%%%%%%%%%%%%%%%%%%%%%%%%%%%%%%
%\subjclass[2010]{65H10, \red{???}}
% 65H10     Nonlinear algebraic or transcendental equations:  Systems of equations 
% 
% 
\keywords{Regeneration, polynomial system solving, homotopy continuation, multidegree}
%%%%%%%%%%%%%%%%%%%%%%%%%%%%%%%%%%%%%%%%%%%%%%%%%%%%%%%%%%%%%%%%%%%%%%%%%%%%%%%%%

\begin{abstract}
We demonstrate our implementation of a 
continuation method as described in \cite{HR2015} for solving polynomials systems. 
Given a sequence of (multi)homogeneous polynomials, the software 
 \texttt{multiregeneration} outputs 
the respective (multi)degree in a wide range of cases and partial multidegree in all others. We use Python for the file processing, while Bertini~\cite{BertiniBook} is needed for the continuation.
Moreover, parallelization options and several strategies for solving structured polynomial systems are available.
\end{abstract} 
%%%%%%%%%%%%%%%%%%%%%%%%%%%%%%%%%%%%%%%%%%%%%%%%%%%%%%%%%%%%%%%%%%%%%%%%%%%% 
\maketitle 
%%%%%%%%%%%%%%%%%%%%%%%%%%%%%%%%%%%%%%%%%%%%%%%%%%%%%%%%%%%%%%%%%%%%%%%%%%%% 

%%%%%%%%%%%%%%%%%%%%%%%%%%%%%%%%%%%%%%%%%%%%%%%%%%%%%%%%%%%%%%%%%%%%%%%%%%%% 
\vspace{-34pt}
\section{Introduction}\label{s:intro}
Many problems in science, engineering, and mathematics are formulated as solving polynomial systems of equations. 
{Homotopy continuation is a technique from the field of applied 
algebraic geometry which can be used approximate all isolated complex 
solutions to a system.} 
The main idea is to construct a parameterized family of problems, obtain all solutions to a member of this family, and then
{use predictor-corrector methods, e.g. Euler-Newton, to track 
the paths of solutions as the solved problem deforms to the problem we 
want to solve.}
For more  details see 
\cite{AG1990, SW05}.

%One approach to solve these systems is to apply a continuation method from \emph{numerical algebraic geometry}, which is also known as numerical nonlinear algebra. 
%A continuation method, also known as a homotopy method, for  solving a polynomial system takes a polynomial system where the solutions are known, and uses numerical predictor-corrector methods to determine solutions of a target polynomial system.
%For an introduction see~\cite{}.

% There are several types of homotopy methods. 
% We focus on a regeneration method, also called an equation-by-equation method. 
% Other methods such as the polyhedral homotopy method \cite{MR3334764, HS1995,PHC} and 
% monodromy homotopy method \cite{BDLS2018,DHJLLS2019} are also~popular. 
\subsection{Witness set and regeneration background}
A \emph{witness set} is used to analyze algebraic varieties and can be computed using homotopy continuation. 
The witness set of a equidimensional projective variety $X$ in $\mathbb{P}^n$ consists of set of the triple of information $(F,L,W)$
such that
$F$ is a witness system \cite[Def. 1.1]{HR2015} of homogeneous polynomials for $X$,
$L$ is a system of $\dim(X)$ general linear polynomials,
and  $W$ is the set of points $V(L)\cap X$.
{ By Bertini's Theorem, the degree of $X$ is the cardinality of $W$.}
%JIR: I did't understand why this change.
%By Bertini's Theorem, the degree of $X$ is the cardinality of \new{any witness set} $W$.% 

Regeneration \cite{HR2015,HSW2011} is {an 
\emph{intersection approach}} for solving polynomial 
systems, {which uses witness sets to keep track of solutions as 
equations are added one by one.}
%\old{Geometrically, regeneration in a projective space $\mathbb{P}^n$ can be understood as an \emph{intersection approach}.} 
%
Given an equidimensional degree $d$ subvariety $X$ in $\pp^n$ and a degree $m$ hypersurface $\mathcal{H}$ not containing an irreducible component of $X$, 
 by Bezout's Theorem the intersection $X\cap\mathcal{H}$ has a degree bounded from above by  $d\cdot m$.
If $\mathcal{H}$ is defined by a general linear product of $m$ linear forms $r_1,\dots,r_m$ 
and $X$ is positive dimensional, then 
the degree bound $d\cdot m$ is always achieved. 
 
{When $g$ is a homogeneous degree $m$ polynomial, 
for $t \in \mathbb{C}$
the family of polynomials $t \prod_{j=1}^m r_j +(1-t) g$ define
 a family of hypersurfaces of constant degree $m$.} % which we denote by $\mathcal{H}_t$.
%We take  $\mathcal{H}_0:=\calV(g)$ and $\mathcal{H}_\infty:= \calV(  \prod_{j=1}^m r_j )$.
By interpolating between members of this family, we are able to track a witness collection for $X\cap\mathcal{V}(\prod_{j=1}^m r_j)$ to a witness collection for $X\cap\mathcal{V}(g)$ using homotopy continuation.
One pass of regeneration 
computes  a witness set for  $X\cap \mathcal{ V }(g)$
 from a witness set for $X \subset \pp^{n}$ by going through  
 two main steps, as described in Algorithm~$\ref{A:regen-hypersurface}$. 

%%%%%%%%%%%%%%%%%%%%%%%%%%%%%%%%%%%%%%%%%%%%%%%%%%%%%%%%%%%%%%%%%%%%%%%%%%%%%%%%% 
\begin{samepage}
\begin{alg}[Regenerating a projective hypersurface intersection]
  \label{A:regen-hypersurface}
  \mbox{\ }\newline
  {\bf Input:}  A witness set $(F,\, L'\cup \{ \ell \},\, W)$ for the positive equidimensional projective variety $X \subset \pp^{n}$ 
  and
a degree $m$ homogeneous polynomial $g$ defining a hypersurface $\mathcal{V}(g)\subset \mathbb{P}^{n}$ not containing any irreducible component $X$.
  \newline  
  {\bf Output:}  
  A witness set  $(F\cup \{g\}, L', W')$ for  $X \cap \mathcal{ V }(g)$.
  \newline
   {\bf Do:} 
  For $ j \in [ m ]$, 
  track the set of points $W$ to a set of points $W_j$ via the homotopy 
  $\mathbb{C}\times\mathbb{P}^n 
  \ni
    (t,\bfx)\mapsto (F(\bfx),\, L'(\bfx),\, t\ell(\bfx)+(1-t)r_j(\bfx) )
    $, where  $r_{j}\in\mathbb{C}[\bfx]$ is a general linear polynomial. 
    Track the set of points $\cup_{j=1}^m W_j$ to a set of points $W'$
    using the homotopy 
     $\mathbb{C}\times\mathbb{P}^n  
     \ni
 (t,\bfx)\mapsto (F(\bfx),\, L'(\bfx),\,  t \prod_{i=1}^m r_j(\bfx)+(1-t)g(\bfx))
 $   to obtain the witness set  $(F\cup \{ g \},\, L' ,\, W')$ for $X \cap\calV (g)$.
% Maybe remark t goes from 1 to zero. 
\end{alg}
\end{samepage}

%\new{The case of an arbitrary projective variety $X$ and an arbitrary 
%hypersurface $\sV(g)$ can be handled by modifying 
%Algorithm~\ref{A:regen-hypersurface} in the following way. If $Z$ is a codimension 
%$c$ irreducible component of $\mathcal{V}(F)  \cap \mathcal{V}(g)$, then 
%either $Z$ is the intersection of $\sV(G)$ with a codimension $c-1$ 
%irreducible component of $\mathcal{V}(F)$ not contained in $\mathcal{V}(g)$,
%or $Z$ is an irreducible component of $\mathcal{V}(F)$ contained in 
%$\mathcal{V}(g)$. Therefore given a collection of witness sets for the 
%irreducible components of $X$, one can compute a collection of witness 
%sets for the irreducible components of $X \cap \sV(G)$ via a combination 
%of Algorithm~\ref{A:regen-hypersurface} and a membership test 
%\cite[Chapter 15]{SW05} to detect whether a component of $X$ is 
%contained in $\sV(G)$.}

\newpage
Coupled with a membership test \cite[Chapter 15]{SW05} to determine if a witness points is contained in $\mathcal{V}(g)$,
 the previous algorithm is iterated to give a method for computing witness sets of  intersections of hypersurfaces. 
This works because
a codimension $c$ irreducible component of $\mathcal{V}(F)  \cap \mathcal{V}(g)$ 
arises as either 
an intersection of a codimension $c-1$ 
 irreducible of
  $\mathcal{V}(F) $ not contained in $\mathcal{V}(g)$
and the hypersurface $\mathcal{V}(g)$,
or 
an irreducible component of $\mathcal{V}(F) $ contained in~$\mathcal{V}(g)$.

We focus on computing (multi)degree of nonsingular equidimensional varieties $X$, but 
regeneration is also used for the 
(singular and nonsingular) nonequidimensional 
case~\cite[Example 4.12 and Section 5.1]{HR2015}.

\subsection{Regenerating a twisted cubic}
Consider the twisted cubic in $\mathbb{P}^3$ 
given by the 
%$2\times 2$-minors of 
%$\left[\begin{smallmatrix}x_{0} & x_{1} & x_{2}\\
%x_{1} & x_{2} & x_{3}
%\end{smallmatrix}\right]$. 
polynomials 
%\[\langle -x_2^2+x_1*x_3, -x_1*x_2+x_3*x_0, -x_1^2+x_2*x_0\rangle.\]
$(x_1 x_3 -x_2^2, x_3 x_0  -x_1 x_2 , x_2 x_0 -x_1^2).$
The figures below illustrate the regeneration procedure. 
Regeneration computes a witness set for a quadratic surface,
 a reducible quartic curve, and twisted cubic respectively.  
The last intersection is determined by evaluations of the last polynomial at the witness points rather than continuation. 
See \cite[Twisted cubic]{Tutorial} for more details on this example.

\begin{figure}[htb!]
\centering
\includegraphics[scale=.11]{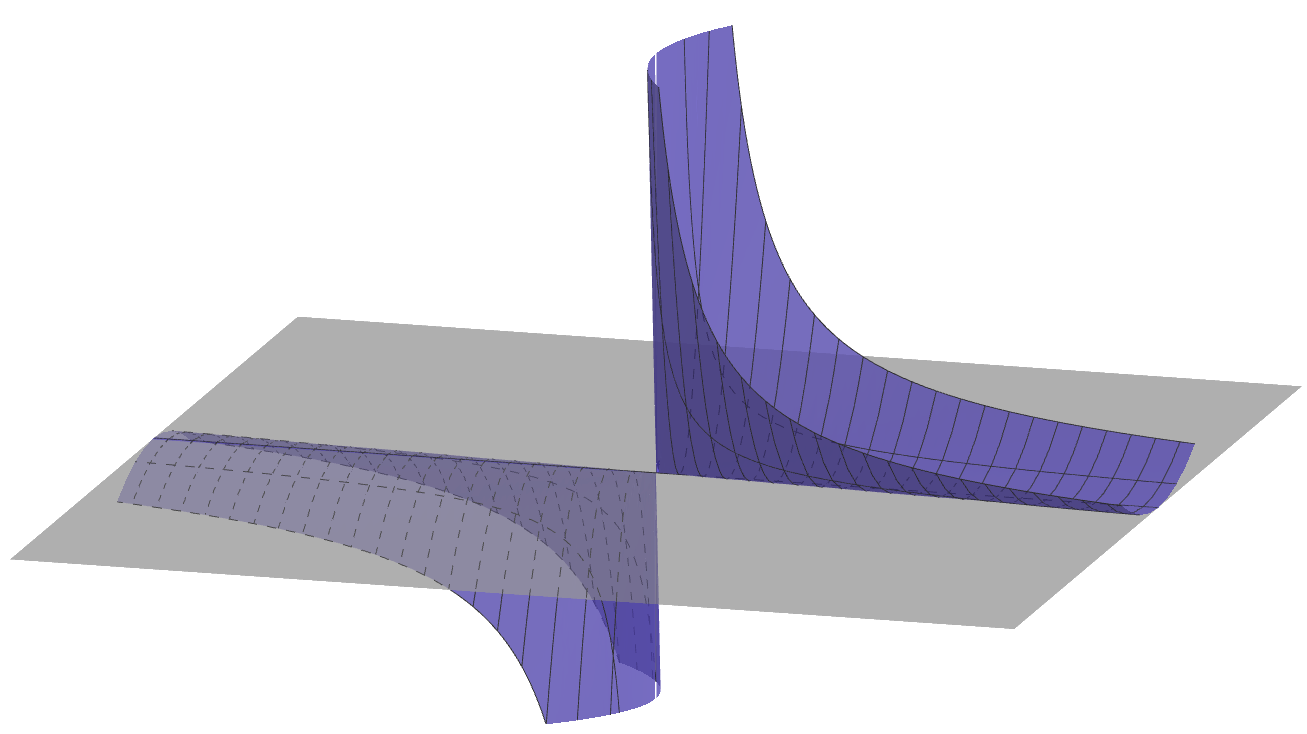}
\includegraphics[scale=.11]{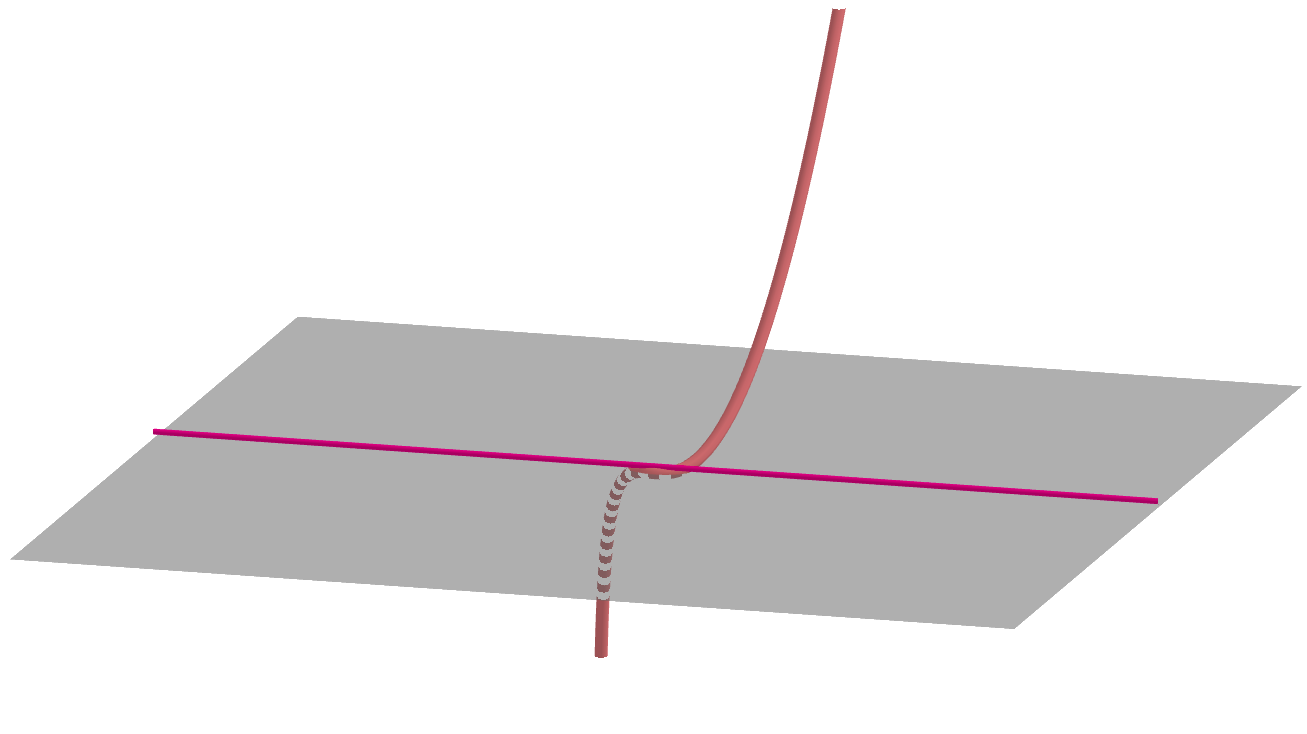}
\includegraphics[scale=.11]{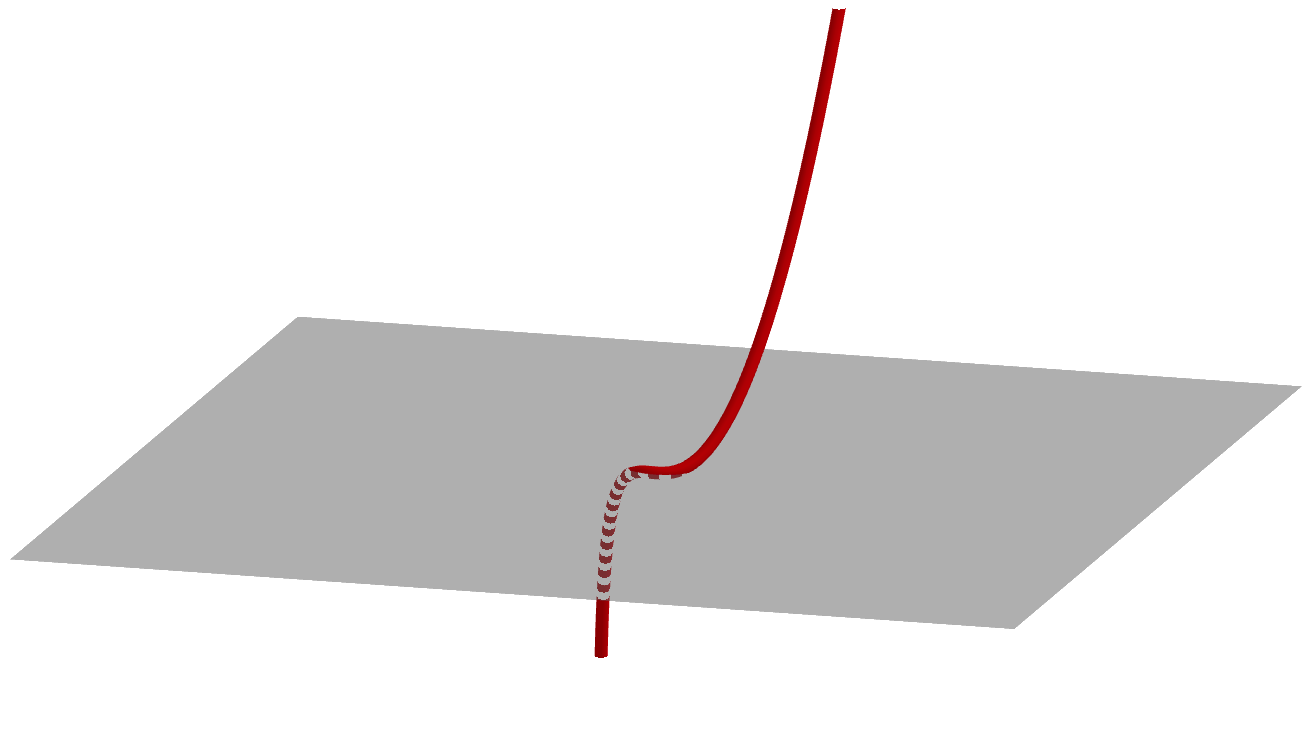}
\end{figure}

Like how regeneration is used to compute the degree of a projective variety, 
\emph{multiregeneration} is used to compute the \emph{multidegree}~\cite{MS05} of a subvariety $X$ of
$\mathbb{P}^{n_1}\times \cdots \times \mathbb{P}^{n_k}$. 
For our purposes it suffices to think of multidegree as a map 
%from the integer lattice to $\mathbb{N}_{\greq}$
$\mathbb{N}^k_{\geq 0}\to\mathbb{N}_{\geq0}$
 that takes an integer vector $\bfe$ with entries summing to $\dim(X)$ to the number of points of intersection of a certain linear space $V(L_\bfe)$ with $X$ and zero otherwise. 
Specifically, $L_\bfe$ is a set of $e_1$ general linear polynomials in $\bfx_1$,
 $e_2$ general linear polynomials in $\bfx_2$, and so on, where $\bfx_i$ are the coordinates of $\mathbb{P}^{n_i}$.

%Just as the degree of a projective variety is determined by intersecting with a general  line space, 
%multidegrees of a subvariety of a product of projective spaces is determined 
%by various intersections with products of general linear spaces. 

\subsection{Regeneration in $\mathbb{P}^3\times\mathbb{P}^1$ example}
Consider the curve $X$ in $\mathbb{P}^3\times\mathbb{P}^1$ given by the ideal
 $\langle x_0y_0 + x_1y_1, x_1y_0 + x_2y_1, x_2y_0 + x_3y_1 \rangle$.
%\ColinDoes{compute this example with equations $F$ and homvariable groups 
%$x_0,..,x_3$ and $y_0,y_1$,
%with our software ( < 1 hr)}
The multidegree of $X$ is determined in Macaulay2~\cite{M2} by the function $\texttt{multidegree}$  using symbolic computation.
%With numerical methods the saturation is trivial because it is carried out by evaluation of a polynomial.
% Too long. 
%\begin{verbatim}
%i1 : R = QQ[x_0..x_3]**QQ[y_0,y_1];
%i2 : I = ideal(x_0,x_1,x_2,x_3)*ideal(y_0,y_1); -- Irrelevant ideal
%i3 : F = ideal(x_0*y_0 + x_1*y_1, x_1*y_0 + x_2*y_1, x_2*y_0 + x_3*y_1);
%i4 : F' =ideal(0_R);
%i5 : apply(F_*,g -> ( F' = saturate(F' + ideal(g), I); multidegree F'  ))      
%        3     2     3     2     3     2
%o27 = {T  + 3T T , T  + 3T T , T  + 3T T }
%        0     0 1   0     0 1   0     0 1
%\end{verbatim}

%Short version
%restart
%R=QQ[x_0,x_1,x_2,x_3]**QQ[y_0,y_1]
%X =ideal(matrix{{y_0,y_1}}*matrix{{x_0,x_1,x_2},{x_1,x_2,x_3}})
%I = ideal (y_0,y_1)*ideal(x_0,x_1,x_2,x_3)
%multidegree saturate( X,I)
%
%ideal(matrix{{y0,y1}}*matrix{{x0,x1,x2},{x1,x2,x3}})
%
%multidegree saturate(X,I)
%restart
%R = QQ[x_0..x_3]**QQ[y_0,y_1];
%I = ideal(x_0..x_3)*ideal(y_0,y_1)
%multidegree saturate(ideal(x_0*y_0 + x_1*y_1, x_1*y_0 + x_2*y_1, x_2*y_0 + x_3*y_1),I)
%%i1 : R = QQ[x_0..x_3]**QQ[y_0..y_1];
%minors(2,matrix{x0,x1,x2})
\begin{leftbar}
\begin{verbatim}
i1 : (n1,n2)=(3,1); R = QQ[x_0..x_n1]**QQ[y_0..y_n2];
i2 : I = ideal(x_0,x_1,x_2,x_3)*ideal(y_0,y_1); -- Irrelevant ideal
i3 : F = ideal(x_0*y_0 + x_1*y_1, x_1*y_0 + x_2*y_1, x_2*y_0 + x_3*y_1);
i4 : multidegree saturate(F,I)      
      3     2
o4 = T  + 3T T
      0     0 1
\end{verbatim}
\end{leftbar}

With our software \texttt{multiregeneration}, each coefficient of the multidegree is exhibited by a witness set. This  is called a witness collection for the multidegree in \cite{Toolkit}. Our implementation computes each of these witness sets.
% representatives. 
 The input is described in the next section, but the displayed output 
%(with the line indexing the variable groups removed to save space) 
(with some formatting) 
is reproduced below.  
 The rows containing numbers
are in bijection with
the terms of \texttt{o4}. Namely, 
the row containing numbers $(c,e_1,e_2)$ is taken to the term $\mathtt{c\frac{T_0^{n_1}T_1^{n_2}}{T_0^{e_1}T_1^{e_2}}}$.
\begin{leftbar}
\begin{verbatim}
python multiregeneration.py
| # smooth isolated solutions  | # of general linear equations |
| found                        | added with variables in group |
  3                              1  0  
  1                              0  1  
\end{verbatim}
\end{leftbar}

%\begin{remark}
\section{Implementation features and examples}
\subsection{High-Throughput Computing Paradigm}
It is natural to interpret regeneration as the exploration of a graph. 
Just as one can search a graph depth first or breadth first, so can one perform multiregeneration. 
The depth first approach prioritizes finding solutions to the system over solving equation by equation in a breadth first approach. 
Previous implementations have taken a breadth first approach, while to our knowledge, we are the first to allow for a depth first approach, which 
is useful 
%in applications to 
to get immediate partial feedback on the original~system. 
 % cat dimension_summary.txt                                   

% Screenshot instead of verbatim.
% \begin{figure}[htb!]
% \centering
%     {\includegraphics[scale=.33]{output.png}}
%   \vspace{-30pt}% There was a lot of White space so I moved the caption up.
%   \caption{Output and file structure}\label{fig:output}
%\end{figure}

The implementation is written in {Python} and design choices have been made to make use of \emph{High-Throughput Computing} (HTC) over 
\emph{High-Performance Computing} (HPC). 
HTC is ideal for jobs that do not communicate with each other 
while HPC is better when rapid communication is needed to perform the computations. 
For this reason, we only focus on regular solutions to avoid the communication needed in comparing endpoints of paths for handling higher multiplicity. 
%A solution is identified by a random tag (as illustrative above) to avoid the need of any sorting. 
Moreover, HTC makes use of many computing resources over long periods of time whereas HPC make intense use of large amounts of computing resources in relatively short time periods. 
Since every endpoint of a completed homotopy is saved as a file with each line giving a~coordinate, we are able to access intermediate results even if the software is running for months. 
For example, the command \texttt{tree} prints the current status of completed solutions.

\begin{leftbar}
{\scriptsize
\begin{verbatim}
tree run/_completed_smooth_solutions   # Output is slightly modified for illustrative purposes. 
|------ depth_0                        # Terminal tree command displaying software's results as a graph.
|   |----- depth_0_gens_1_dim_2_1_varGroup_1_regenLinear_1_pointId_429439718721_285170369818
|   |----- depth_0_gens_1_dim_3_0_varGroup_1_regenLinear_1_pointId_429439718721_494593912469
|----- depth_1
|   |----- depth_1_gens_1_1_dim_1_1_varGroup_1_regenLinear_1_pointId_285170369818_258141170677
|   |----- depth_1_gens_1_1_dim_2_0_varGroup_1_regenLinear_1_pointId_285170369818_257010786796
|   |----- depth_1_gens_1_1_dim_2_0_varGroup_1_regenLinear_1_pointId_494593912469_916362171011
|----- depth_2
  |----- depth_2_gens_1_1_1_dim_0_1_varGroup_1_regenLinear_1_pointId_258141170677_929159838948
  |----- depth_2_gens_1_1_1_dim_1_0_varGroup_1_regenLinear_1_pointId_257010786796_156506717710
  |----- depth_2_gens_1_1_1_dim_1_0_varGroup_1_regenLinear_1_pointId_258141170677_604647130850
  |----- depth_2_gens_1_1_1_dim_1_0_varGroup_1_regenLinear_1_pointId_916362171011_957285449047
\end{verbatim}
}
\end{leftbar}

\begin{remark}
Bertini has its own regeneration option for subvarieties  of 
$\mathbb{P}^n$ and $\mathbb{C}^n$ in addition  to  zero-dimensional subvarieties of a product of projective or affine spaces. 
\end{remark}

\subsection{Input files and features}
The software \texttt{multiregeneration} needs four input files to run. 
The first input file, \texttt{inputFile.py}, contains the (multi)degree information about each polynomial of the system to be solved. 
Furthermore, the user can set options 
%\texttt{TorusOnly}, 
%\texttt{NonZeroCoordinates},
to access additional features, such as %implementing side conditions, 
restricting a witness set to a hypersurface, 
working with the algebraic torus,
and returning specific coefficients of the multidegree. 
\begin{leftbar}
\begin{verbatim}
## The text in inputFile.py to run the twisted cubic example.
degrees = [[2], [2], [2]]		# The degrees of the three polynomials.
verbose = 1		# Change to zero to display nothing.
algebraicTorusVariableGroups = [0]	# List of variable groups where
# returned solutions have nonzero coordinates
maxProcesses = 1 # Change to N to use N processes in parallel
\end{verbatim}
\end{leftbar}
%?????Processes =  1 # Change to N to use N processes in parallel.

%The variable "degrees" must be initialized to a list of lists, where the 'th element of the 'th list is the degree of the 'th function in the th "variable group." 
%Bertini uses "I" to denote . The use of this symbol as a variable is not allowed. To specify the value of a constant, say is , put in a single line in the "bertiniInput\_equations" file "c = 2.2".

The remaining files needed to run the software are used to automate calling Bertini during the multiregeneration. 
The files are (1)  \texttt{bertiniInput\_variables} to describe the variable groups used (also allows for the affine case), 
(2)  \texttt{bertiniInput\_equations}  to describe the system of polynomials to be solved,
and  (3)  \texttt{bertiniInput\_trackingOptions} to customize the continuation tracker  (see Appendix A in the Bertini user manual).
%These files are written in the $C$-like syntax used by the Bertini software.
%. In the "bertiniInput\_variables" file, the unknowns of our system of polynomials are . In thefile, our system of polynomials is described by a line beginning with "function" to set the polynomials whose common roots we aim to describe followed by one equation per line to define the polynomials in an expression of the unknowns.

\begin{leftbar}
\begin{verbatim}
# The text in bertiniInput_variables for the twisted cubic example.
# Use variable_group to work with affine coordinates and 
#  multiple lines to work with several variable groups. 
hom_variable_group x_0, x_1, x_2, x_3;

# The text in bertiniInput_equations. 
# It follows the same syntax as Bertini.
function f1, f2, f3;
f1 = x_1^2 - x_0*x_2; f2 = x_2^2 - x_1*x_3; f3 = x_0*x_3 - x_1*x_2;

# The text in bertiniInput_trackingOptions is to achieve 
#  tighter convergence tolerance for the endgame. 
FinalTol: 1e-12;
\end{verbatim}
\end{leftbar}

For additional illustrative examples see~\cite{Tutorial}.
These cover a range of topics: 
economics (Nash equilibria~\cite{MR2551729}),
machine learning  (deep linear networks~\cite{Tingting}), 
and
statistics (maximum likelihood estimation~\cite{HRS14}).
%kinematics (Alt's problem)\cite{},
%and computer~vision.
\subsection*{Acknowledgments}
Our work was partially supported by the National Science Foundation through grant number DMS-1502553.
Research of Rodriguez supported in part by the College of Letters and Science.The authors would like to thank Melanie Wood and Nigel Boston for their support of the CURL 2019 Summer program. 
The authors would like to thank Steve Goldstein, Diego Cifuentes and Tingting Tang for helpful comments.

\bibliographystyle{abbrv} % abbreviates first and middle names
\bibliography{regen}
%\pagebreak

%\ColinDoes{ translate the file from here to the preferred format syle by ISSAC. (<.5 hr)}
%
%
%\ColinDoes{state exactly what our implementation does in the Github readme. 
%This should clearly describe what this sentence is trying to describe: 
%``We implement a version in \cite{HR2015} for complete intersections where the hypersurface $\cH_m$ intersects the variety $\cap_{i<m}\cH_i$ transversally for each irreducible component. "}
%
%
%\ColinDoes{on June 2, check every reference and citation. 
%This means looking it up online through the UW Madison library website to get access. 
%Bibtex should be done with Mathscinet whenever possible. 
%Check that the Theorems numbers and example numbers agree with the published versions.  (< 2 hrs)}
%
\end{document}